\documentclass[10pt,conference]{ieeeconf}
\IEEEoverridecommandlockouts
\usepackage{cite}
\usepackage{amsmath,amssymb,amsfonts,mathtools}
\usepackage{algorithmic}
\usepackage{graphicx}
\usepackage[export]{adjustbox}
\usepackage{tikz}
\usepackage{textcomp}
\usepackage{xcolor}
\usepackage[font=small]{caption}
\def\BibTeX{{\rm B\kern-.05em{\sc i\kern-.025em b}\kern-.08em
		T\kern-.1667em\lower.7ex\hbox{E}\kern-.125emX}}

\usepackage{amsthm}

\renewcommand{\(}{\left(}
\renewcommand{\)}{\right)}
\renewcommand{\[}{\left[}
\renewcommand{\]}{\right]}

\newcommand{\reals}{\mathbb{R}}
\newcommand{\transpose}{^\textrm{T}}
\newcommand{\continuous}{\mathcal{C}}
\DeclareMathOperator*{\argmin}{arg\,min}
\DeclareMathOperator*{\argmax}{arg\,max}

\newtheorem{theorem}{}

\newtheorem{assumption}{}

\newtheorem{lemma}{}

\newtheorem{definition}{}

\newtheorem{remark}{}

\newtheorem{problem}{}

\newtheorem{example}{}

\usepackage{cancel}
\usepackage[normalem]{ulem} 
\usepackage{verbatim}
\usepackage{url}
\usepackage{hyperref}
\usepackage{placeins}
\usepackage{soul}

\definecolor{lcolor}{rgb}{0,0,0.6}
\definecolor{Blue}{rgb}{0,0,0.8}
\hypersetup{
	colorlinks=true,
	linkcolor=lcolor,
	filecolor=magenta,      
	urlcolor=cyan,
}

\begin{document}

\title{{\fontsize{15.9}{16} \bf High Relative Degree Control Barrier Functions Under Input Constraints}}%

\author{Joseph Breeden \and Dimitra Panagou\thanks{The authors are with the Department of Aerospace Engineering, University of Michigan, Ann Arbor, MI, USA. Email: \texttt{\{jbreeden,dpanagou\}@umich.edu}.}\thanks{The authors thank the National Science Foundation Graduate Research Fellowship Program for supporting this research.}}

\maketitle

\begin{abstract}
	
This paper presents methodologies for ensuring forward invariance of sublevel sets of constraint functions with high-relative-degree with respect to the system dynamics 
and in the presence of input constraints. We show that such constraint functions can be converted into special Zeroing Control Barrier Functions (ZCBFs), which, by construction, generate sufficient conditions for rendering the state always inside a sublevel set of the constraint function in the presence of input constraints. 
We present a general form for one such ZCBF, as well as a special case applicable to a specific class of systems.
We conclude with a comparison of system trajectories under the two ZCBFs developed and prior literature, and a case study for an asteroid observation problem using quadratic-program based controllers to enforce the ZCBF condition.
	
\end{abstract}

\section{Introduction}

Control Barrier Functions (CBFs) have recently gained popularity across disciplines for control synthesis in safety-critical systems. However, the problem of designing CBFs for high-relative-degree ($r \geq 2$) constraint functions under control input constraints remains an open question except for specific systems \cite{First_CBF,Grizzle_SeminalPaper,Robot_D2CBFs,Integrator_CBFs}. 
In this paper, we seek to generalize these results, and in particular the methods in \cite{Integrator_CBFs,My_Paper,backup_controller}, to a wider class of systems. 

System safety is often formulated as a set invariance problem, wherein the controller must render the state trajectories always inside a designated \textit{safe set}. In this work, the safe set is specified as the zero sublevel set of some \textit{constraint function}. If the constraint function is of relative-degree $r=1$ with respect to the system dynamics, then the set can be rendered forward invariant if the constraint function is also a CBF \cite{Grizzle_SeminalPaper}. For control-affine dynamics, this leads to an affine condition on the control input, which can be applied pointwise to yield either explicit control laws \cite{First_CBF}, or optimization-based control laws \cite{Grizzle_SeminalPaper,Backstepping_CBF,Exponential_CBF,CBFs_ManyClassK} that render the safe set forward invariant. However, if the constraint function is of relative-degree $r \geq 2$ and the uncontrolled dynamics allow trajectories to leave the safe set, then the constraint function alone cannot be a CBF. To recover such a control-affine condition, we seek to construct CBFs that are composed of the constraint function and its derivatives, and whose sublevel sets are subsets of the safe set.

Several papers develop methods to convert high-relative-degree constraint functions ($r \geq 2$) into CBFs, including using compositions with bounded monotonic functions \cite{Grizzle_SeminalPaper}, backstepping \cite{Backstepping_CBF}, feedback linearization and pole placement \cite{Exponential_CBF,XU2018195}, or by defining safe sets for every order of derivative of the constraint function \cite{CBFs_ManyClassK}. The approach in \cite{My_Paper} bypasses the creation of a new CBF, but develops a condition on the first controllable derivative that fills the same role as the conditions in \cite{Backstepping_CBF,Exponential_CBF}. All these approaches lead to conditions on the control input that are potentially infeasible if the set of valid control inputs is bounded. In practice, input constraints may be satisfied within the frameworks of \cite{Exponential_CBF,Backstepping_CBF,XU2018195,Grizzle_SeminalPaper} for certain trajectories by tuning (e.g. choosing different poles using the method in \cite{Exponential_CBF}), but these approaches are only provably feasible (and hence, provably safe) if the control set is $\reals^m$. The work in \cite{CBFs_ManyClassK} improves upon this by defining a subset of the safe set that is controlled forward invariant in the presence of input constraints. However, choosing appropriate class-$\mathcal{K}$ functions to satisfy the feasibility requirements in \cite[Def. 7]{CBFs_ManyClassK} may not be straightforward. 

Conditions for safety under input constraints for the $n$-integrator system are introduced in \cite{Integrator_CBFs}. For certain other systems, a second CBF that guarantees satisfaction of control input constraints for all future times can also be introduced \cite{Grizzle_SeminalPaper,Robot_D2CBFs}. For more general systems, \cite{backup_controller,backup_implementation} recently developed an approach (not specific to high-relative-degree) wherein a small known \textit{backup set} is expanded to the set of states which can reach the backup set in a finite time horizon under input constraints. The work in \cite{Performance_Func} is similar to \cite{backup_controller,backup_implementation}, but generalizes the approach to infinite time horizon.

This paper addresses the problem of designing CBFs for high-relative-degree constraint functions with guaranteed safety under input constraints for a broader class of systems than in \cite{First_CBF,Integrator_CBFs,Grizzle_SeminalPaper,Robot_D2CBFs,backup_controller,backup_implementation} using extensions to the approach in \cite{Performance_Func}. We introduce two strategies for generating CBFs that by construction can be rendered nonpositive in the presence of input constraints. The first strategy uses a predefined nominal control law similar to \cite{backup_controller,backup_implementation,Performance_Func}, but does not require a predefined backup set (as in \cite{backup_controller,backup_implementation}), and does not require the searched time horizon to contain a unique maximizer of the constraint function (as in \cite{Performance_Func}). The second strategy simplifies the first and generalizes \cite{Integrator_CBFs} to any system for which there exists a minimum control authority over the $r$th derivative of the constraint function everywhere in the safe set. Similar to \cite{CBFs_ManyClassK}, the resultant controlled forward invariant sets are subsets of the safe set, and could be smaller or larger than the corresponding sets obtained in \cite{CBFs_ManyClassK}. However, unlike \cite{CBFs_ManyClassK}, these methods by construction respect input constraints without tuning other parameters of the CBF. We then simulate the new strategies on an obstacle avoidance problem for a double-integrator system, and demonstrate the second strategy on a safety-critical spacecraft control problem for asteroid observation.


\section{Preliminaries} \label{sec:preliminaries}

\subsection{Notation}

Let $\continuous^r$ be the set of $r$-times continuously differentiable functions. 
Let $\emptyset$ denote the empty set and $\partial S$ the boundary of set $S$. Let $|| \cdot ||$ be the $2$-norm, and $||\cdot||_\infty$ be the $\infty$-norm. For $v\in\reals^n$, let $v^\perp = \{w\in \reals^n : v\transpose w = 0\}$. A function $\alpha : \reals \rightarrow \reals$ belongs to extended class-$\mathcal K$, denoted $\alpha \in \mathcal{K}$, if it is strictly increasing and $\alpha(0) = 0$. Let $L_f h (x)= \frac{\partial h}{\partial x} f(x)$ denote the Lie-derivative of a function $h : \reals^n \rightarrow \reals$ with respect to a function $f : \reals^n \rightarrow \reals^{n\times m}$ at the point $x$. Let $L_f^2h(x) = L_f(L_f h(x))$ and $L_f^r h(x) = L_f (L_f^{r-1} h(x))$. 

\subsection{Model and Problem Formulation}

We consider the control-affine system
\begin{equation}
	\dot{x}(t) = f(x(t)) + g(x(t))u(t) \,, \label{eq:model}
\end{equation}
with state $x\in\reals^n$, control input $u\in U \subseteq \reals^m$ where $U$ is compact, and functions $f:\reals^n\rightarrow \reals^n, g:\reals^n\rightarrow\reals^{n\times m}$ where $f,g \in \continuous^r$.


\begin{definition}
A function $h : \reals^n \rightarrow \reals$ is said to be of \textnormal{relative-degree $r$} with respect to the dynamics \eqref{eq:model} if
\begin{enumerate}
	\item $h \in \continuous ^r$,
	\item $L_g L_f^k h(x) = 0, \; \forall x\in\reals^n, \; \forall k=0,1,\cdots r-2$, and
	\item $\exists C \subseteq \reals^n, C \neq \emptyset$ such that $L_g L_f^{r-1} h(x) \neq 0,\forall x \in C$.
\end{enumerate}
\end{definition}
\newcommand{\highdegree}{\mathcal{G}}

We denote the space of all functions of relative-degree $r$ with respect to a given system as $\highdegree^r$. Let $h : \reals^n \rightarrow \reals, h\in \highdegree^r$, and define a safe set $S$ as
\begin{equation}
	S \triangleq \{x\in\reals^n \mid h(x) \leq 0\} \,. \label{eq:safe_set}
\end{equation}
We call $h$ the \textit{constraint function} for set $S$. If $h \in \highdegree^1$, then a sufficient condition for forward invariance of $S$ is that $h$ is a Zeroing Control Barrier Function, defined as follows.
\begin{definition}[{\hspace{-0.2pt}\cite[Def. 5]{Grizzle_SeminalPaper}}] 
	\label{def:zcbf}
	For the system \eqref{eq:model}, a continuously differentiable function $h:\reals^n\rightarrow \reals$ is a \textnormal{zeroing control barrier function (ZCBF)} on set $S$ in \eqref{eq:safe_set} if $\exists \alpha\in\mathcal{K}$ such that
	\begin{equation}
		\inf_{u \in U} \[ L_f h(x) + L_g h(x) u - \alpha (-h(x))\] \leq 0 ,\; \forall x \in S\,. \label{eq:cbf}
	\end{equation}
\end{definition}
\noindent Forward invariance of $S$ is then guaranteed as follows.
\begin{lemma}[{\hspace{-0.2pt}\cite[Cor. 2]{Grizzle_SeminalPaper}}] \label{prior:cbf_invariance}
	If $h$ is a ZCBF for $S$ under \eqref{eq:model}, then any Lipschitz continuous controller $u\in U$ such that
	\begin{equation}
		\dot{h}(x,u) = L_f h(x) + L_g h(x) u(x) \leq \alpha (-h(x)), \;\forall x \in S \label{eq:cbf_condition}
	\end{equation}
	 will render the set $S$ forward invariant.
\end{lemma}

\noindent Note that \ref{def:zcbf} can be relaxed if $h$ is only differentiable \cite[Thm. 2.1]{Additive_CBFs}, and \ref{prior:cbf_invariance} can be extended to non-Lipschitz controllers if \eqref{eq:model} admits a unique solution \cite{BLANCHINI19991747}.

If $h \in \highdegree^r$ for $r \geq 2$, then $L_g h(x) u \equiv 0$, so $h$ must satisfy \eqref{eq:cbf_condition} for all $u \in U$ to be a ZCBF. This is not useful for control design. Thus, the objective of this paper is as follows.

\begin{problem} \label{problem}
	Given a constraint function $h \in \highdegree^r$ for $r \geq 2$, such that $\exists x_0 \in \partial S : L_f h(x_0) > 0$, and a compact control set $U$, develop functions $H : \reals^n \rightarrow \reals$ such that $S_H = \{x\in\reals^n \mid H(x) \leq 0\}$ is a subset of $S$, and $H$ is a ZCBF on $S_H$.
\end{problem}

A solution $H$ to \ref{problem} defines a subset $S_H$ of the safe set that can be rendered forward invariant. \ref{problem} can be addressed by methods such as those in \cite{CBFs_ManyClassK} with a proper selection of $\alpha_i\in\mathcal{K}, i=1,2,\cdots,r$ (or without input constraints for any $\alpha_i\in\mathcal{K}$). In contrast, the approaches in this paper always satisfy the input constraints.

Lastly, similar to \cite{Integrator_CBFs,Performance_Func}, we define the flow operator $\psi_h(t; x, u)$ for $t \geq 0$ as the value $h(y(t))$ resulting from the initial value problem $\dot{y} = f(y)+g(y)u, y(0) = x$ under the control law $u$. 
Also, let $\psi_x(t;x,u)$ denote the value of the state $y(t)$ according to the same initial value problem. 


\section{Methodologies}
\label{sec:main_result}

\subsection{General Case} \label{sec:general}

For $h\in \highdegree^r$, we refer collectively to the derivatives $\dot{h},\ddot{h}, \cdots h^{(r-1)}$ which are not explicit functions of $u$ as \textit{generalized inertia}, by analogy to inertia in kinematic systems. 
To ensure safety when $r \geq 2$, a controller must be able to dissipate this generalized inertia, i.e. ensure $h^{(k)}(x(t)) \leq 0, \forall k\leq r$ for all $t$ such that $h(x(t)) = 0$.

Our approach is to examine the system response forward in time according to its generalized inertia.
Suppose that $h(x) < 0$ and $h^{(k)}(x) > 0$ for one or more $k \in \{1,2,\cdots r\}$. Then we seek to determine how large each $h^{(k)}$ can be allowed to grow before there is no allowable control input under which the trajectory stays within the safe set at some future time. However, analyzing all possible future trajectories (i.e. all possible control laws) for safety is intractable, so instead suppose that we have a predefined control law $u^* : \reals^n \rightarrow U$ that attempts to drive the state towards the interior of $S$ (called a ``nominal evading maneuver'' in \cite{Performance_Func}). For example, if $U$ is a closed ball, one might choose
\\ \vspace{-8pt}
\begin{equation}
	u^*_\textrm{ball}(x) \triangleq \argmin_{u \in U} L_gL_f^{r-1} h(x) u \,,  \label{eq:opt_control_law}
\end{equation}
\vspace{-10pt} \\
which pointwise minimizes $h^{(r)}$.

\vspace{-3pt}\begin{assumption} \label{as:uniqueness}
	For any $x(0)\in S$, the system \eqref{eq:model} admits a unique solution $x(t), t\geq 0$ under the control law $u^*(x)$.
\end{assumption}\vspace{-4pt}%

Denote the evolution of $h$ under $u^*$ from initial condition $x(0)$ as $\psi_h(t; x(0), u^*)$. Then there exists at least one safe and feasible trajectory from $x(0)$ if $\psi_h(t, x(0), u^*) \leq 0, \forall t \geq 0$. Thus, we introduce a new function,
\\ \vspace{-8pt} 
\begin{equation}
    H (x) \triangleq \sup_{t\geq 0} \psi_h (t; x, u^*) \,, \label{eq:new_cbf}
\end{equation}
\vspace{-10pt} \\
which we seek to render nonpositive (see also \cite[Eq. 14]{Performance_Func}). 

\vspace{-1pt}\begin{assumption} \label{as:existence}
	For all $x\in S$, $H(x)$ exists and is finite and differentiable.
\end{assumption}\vspace{-2pt}%
\noindent That is, we assume that the trajectories under $u^*$ dissipate the generalized inertia so that $\psi_h$ is upper bounded, and this upper bound is regular with respect to $x$.

\begin{figure}
	\centering
	\includegraphics[width=\columnwidth]{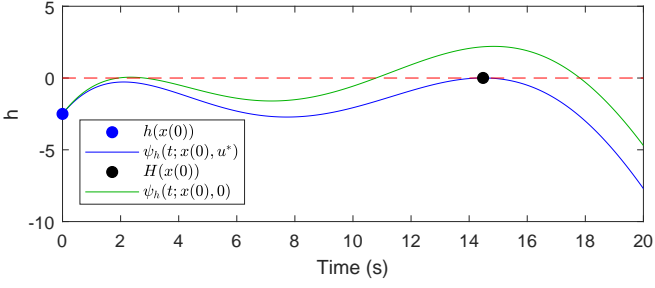}
	\caption{A potential trajectory of $\psi_h(t; x, u^*)$ under control input $u=u^*(x)$, and unforced evolution $\psi_h(t; x, 0)$ under $u=0$ for comparison. We are interested in whether the point $H(x)$ exceeds 0, which would imply that under our proposed controller $u^*$, there exists $t$ such that $h(x(t)) > 0$. \vspace{-11pt}}
	\label{fig:visual}
\end{figure}

The relationship between $h(x)$, $\psi_h(t; x, u^*)$, and $H(x)$ is visualized in Fig.~\ref{fig:visual}. 
In this example, $h(x(0)) < 0$, but there exist $t > 0$ such that $h(x(t)) > h(x(0))$ under the control law $u^*$. We are interested in the value of $H(x(0))$ (black dot in Fig.~{\ref{fig:visual}}), because $H(x(0)) > 0$ would imply that there exists $t>0$ such that $h(x(t)) > 0$ under $u^*$, which means $u^*$ does not render the system safe from $x(0)$.

If the supremum in \eqref{eq:new_cbf} is achieved for finite $t$, denote the set of time instances that maximize $\psi_h$ as
\begin{equation}
t_c(x) \triangleq \Big\lbrace \argmax_{t\geq 0} \psi_h(t; x, u^*) \Big\rbrace \,, \label{eq:def_tc}
\end{equation}
which allows us to equivalently write
\begin{equation}
H(x) = \psi_h(t_{c,0}; x, u^*), \; \forall t_{c,0} \in t_c(x) \,.
\end{equation}
Note that the set $t_c(x)$ may contain more than one element. 
One can see this visually for the example in Fig.~\ref{fig:visual}, for which a slight perturbation 
 might result in the blue line having two maximizers. 
Finally, define the set rendered safe by $u^*$ as
\begin{equation}
	S_H \triangleq \{x\in \reals^n \mid H(x) \leq 0 \} \,. \label{eq:SH_set}
\end{equation}
By definition, $H(x) \geq h(x), \forall x \in \reals^n$, so $S_H \subseteq S$.


For brevity, in this section let $\psi_h(t) = \psi_h(t; x, u^*)$ and $\psi_x(t) = \psi_x(t; x, u^*)$. We are now ready to state the first main result of this paper.

\begin{theorem} \label{thm:most_general}
	The function $H$ in \eqref{eq:new_cbf} is a ZCBF on the set $S_H$ in \eqref{eq:SH_set} for the control set $U$, provided $S_H \neq \emptyset$.
		
	\begin{proof}
		$H$ is a ZCBF if it meets the ZCBF condition \eqref{eq:cbf}. 
		We consider three cases depending on the elements of $t_c(x)$. 
		
		First, consider if $0\in t_c(x)$. Then it must hold that $\dot{h}(x) \leq 0$ (where $\dot{h}$ is independent of $u$ if $h\in\highdegree^r$ for $r\geq 2$) for $0$ to be a maximizer of $\psi_h(t)$. By definition, $\psi_h(0) = h(x)$ and thus $L_f \psi_h (0) + L_g \psi_h(0) u = L_f h(x) + L_g h(x) u = \dot{h}(x) \leq 0$. Thus, the ZCBF condition is satisfied for any $u$ if $h \in \highdegree^r, r\geq 2$, and for the input $u^*$ if $h\in\highdegree^1$.
		
		Next, consider if $t_{c,0}\in t_c(x)$ where $t_{c,0} > 0$. Then $t_{c,0}$ is a maximizer on an open interval, so a necessary condition is $\frac{\partial \psi_h}{\partial t}\big |_{t_{c,0}} = 0$. Define $\kappa \triangleq \psi_h(t_{c,0}; x, u^*)$. It follows that
		\begin{equation}
			\kappa = \psi_h(t_{c,0}; x,u^*) = \psi_h(t_{c,0}-\tau ; \psi_x(\tau; x, u^*), u^*) \,,
		\end{equation}
		$\forall \tau \hspace{-2pt}\in\hspace{-2pt} [0,t_{c,0}]$. 
		Since $\kappa$ is constant w.r.t $\tau$, its derivative satisfies
		\begin{multline}
			\frac{d}{d\tau}[\kappa] = \frac{\partial \psi_h}{\partial t} \frac{\partial (t_{c,0}-\tau)}{\partial \tau} + \frac{\partial \psi_h}{\partial x}\frac{\partial \psi_x(\tau; x, u^*)}{\partial \tau} \implies\\
			0 = - \frac{\partial \psi_h}{\partial t} + \frac{\partial \psi_h}{\partial x}(f(\psi_x(\tau)) + g(\psi_x(\tau)) u^*(\psi_x(\tau)) \,. \label{eq:proof}
		\end{multline}
		Evaluating the above at $\tau = 0$ yields $\frac{\partial \psi_h (t_{c,0}-\tau)}{\partial t}\big|_{\tau=0} = \frac{\partial \psi_h}{\partial t}\big|_{t_{c,0}}  = 0$, and $\psi_x(0) = x$, so \eqref{eq:proof} becomes
		\begin{equation} \label{eq:nonincreasing}
			0 = L_{f(x)} \psi_h(t_{c,0}; x, u^*) + L_{g(x)} \psi_h(t_{c,0}; x, u^*)u^*(x) \,.
		\end{equation}
		Thus, the input $u^*$, which by definition is always in $U$, renders $L_f \psi_h(t_{c,0}) + L_g \psi_h(t_{c,0}) u = 0$, thereby satisfying the ZCBF condition.
		
		Finally, if the supremum in \eqref{eq:new_cbf} is not achieved for finite $t$ (i.e. \eqref{eq:def_tc} does not exist), then choose $\kappa = \lim_{t\rightarrow\infty} \psi_h(t-\tau; \psi_x(\tau;x,u^*), u^*)$. By \ref{as:existence}, the limit exists, so $\lim_{t\rightarrow\infty} \frac{\partial \psi_h}{\partial t} = 0$. Differentiating $\kappa$ w.r.t $\tau$ as in the prior case yields that \eqref{eq:nonincreasing} holds in the limit as $t_{c,0}\rightarrow \infty$, so under the input $u^*$, the ZCBF condition is still satisfied. As an abuse of notation, we denote this case as $\infty \in t_c(x)$.
		
		Since $t_c$ may contain multiple elements, the trajectory of $H$ satisfies $\dot{H}(x,u) \in \lbrace L_f \psi_h(t_{c,0}) + L_g \psi_h(t_{c,0}) u : t_{c,0} \in t_c(x) \rbrace$ (see also \cite[Sec. II]{Additive_CBFs}). $u^*(x)$ is independent of $t_{c,0}$, so the above cases show that there exists a single $u\in U$ ($u=u^*(x)$) that renders every element of this set nonpositive (see also \cite[Thm. 1]{NonSmooth_CBFs}). That means $\forall x \in S_H, \exists u(x) \in U$ s.t. $\dot{H}(x,u(x)) \leq 0 \leq \alpha(-H(x))$ for any $\alpha \in \mathcal{K}$, which is precisely the definition of a ZCBF over $S_H$.
	\end{proof}
\end{theorem}

The immediate consequence of \ref{thm:most_general} is that if $x(0) \in S_H$, then there is a controller such that $x(t)\in S_H, \forall t\geq 0$, and by extension $x(t) \in S, \forall t \geq 0$ since $S_H\subseteq S$. Specifically, \ref{thm:most_general} implies that if $x(0) \in S_H$, then there exists at least one safe trajectory. Unlike in \cite{Performance_Func}, \ref{thm:most_general} also allows for the possibility of $\psi_h$ having multiple maximizers. The following remark then provides a means to calculate $\dot{H}$ for general systems, which we then apply as a condition on the control input using \eqref{eq:cbf_condition}.

\begin{remark} \label{remark:calculate_dphi}
	Suppose that the control law $u^*:\reals^n\rightarrow U$ satisfies $u^*\in\mathcal{C}^1$. For a control input $v\in U$, $\dot{H}$ is given by
	\begin{equation}
		\dot{H}(x,v) = \max_{t_{c,0}\in t_c(x)} \(\frac{\partial \psi_h(t_{c,0})}{\partial x}\(f(x)+g(x)v\)\) \,. \label{eq:Hdot_explicit}
	\end{equation}
	The gradient of $\psi_h(t_{c,0}; x, u^*)$ w.r.t $x$ at a particular $t_{c,0}$ is given by
	\begin{equation}
		\frac{\partial \psi_h(t_{c,0})}{\partial x} = \frac{\partial h(\psi_x(t_{c,0}))}{\partial x}\theta(t_{c,0}) \label{eq:grad_psi}
	\end{equation}
	where matrix $\theta(t)$ is the solution of the initial value problem 
	\begin{equation}
		\begin{aligned}
			\dot{\theta} &= \frac{\partial}{\partial y}[f(y)+g(y)u^*(y)] \big|_{y=z}\theta, \;\;& \theta(0) = I \\
			\dot{z} &= f(z) + g(z)u^*(z) ,\; &z(0) = x
		\end{aligned} \label{eq:the_ode}
	\end{equation}
	where $I$ is the identity matrix.
\end{remark}

\noindent The consequence of \ref{remark:calculate_dphi} is that $u$ must satisfy a condition of the form \eqref{eq:cbf_condition} once for each element of $t_c(x)$. 

While \ref{thm:most_general} theoretically applies to any system for which a nonempty $S_H$ exists, in practice, it may be limited by the requirement to propose a ``good'' $u^*$ (i.e. one which yields a large $S_H$). Also, for most systems, $H(x)$,$\frac{\partial \psi_h}{\partial x}$ will not have explicit expressions. That said, computing $H(x)$,$\frac{\partial \psi_h}{\partial x}$ only requires propagating two ODEs, which can be done efficiently. On the other hand, the advantage of this approach is that if $H(x(0)) \leq 0$, then one immediately knows $S_H$ can be rendered forward invariant under the input constraints. 

Next, we present an alternative approach that avoids the complexity of \eqref{eq:Hdot_explicit}-\eqref{eq:the_ode} but yields a different $S_H$.

\subsection{Special Case of Constant Control Authority} \label{sec:poly}

Instead of allowing for any $u^*$ satisfying Assumptions~1-2, which could make \eqref{eq:new_cbf} difficult to compute, in this section, we set the control input so as to regulate the system to a constant rate of generalized inertia dissipation. That is, choose any $u$ such that $h^{(r)}(x,u)$ is a predefined constant. 
Specifically, we require a control law $u'$ such that \ref{as:uniqueness} holds and $u'(x) \in \boldsymbol\mu(x) \subset U, \forall x \in \reals^n$, where we define:
\begin{equation}
	\boldsymbol\mu(x) \triangleq \big\lbrace u\in U \;\big|\; h^{(r)}(x,u) = -a_{max} \big\rbrace \,, \label{eq:def_uset}
\end{equation}
where $a_{max}$ is a precomputed constant rate of generalized inertia dissipation,
\begin{multline}
    a_{max} \triangleq \max \bigg( \Big\lbrace a\in \reals \;\big|\; \forall x \in S, \exists v \in (L_g L_f^{r-1}h(x))^\perp :
    \\ -\frac{(a + L_f^r h(x))(L_g L_f^{r-1}h(x))}{|| L_g L_f^{r-1} h(x)||^2} + v \in U \Big\rbrace \bigg) \,, \hspace{-6pt} \label{eq:amax}
\end{multline}
assuming an $a_{max} > 0$ exists.
This choice of $u$ is reasonable if $\min_{u\in U} h^{(r)}(x,u)$ does not vary much with $x$, 
but may be overly conservative in other cases, where the system may be able to dissipate generalized inertia at a rate higher than $a_{max}$ except within a small subset of $S$. 

Under any $u \in \boldsymbol{\mu}(x)$, it follows that $h^{(r)}(x,u) = -a_{max}$, so $\psi_h$ has the Taylor expansion:
\begin{equation}
    \psi_h(t; x, u') = \sum_{i=0}^{r-1} \frac{1}{i!}h^{(i)}(x)t^i - \frac{1}{r!}a_{max}t^r \label{eq:phi_poly} \,.
\end{equation}
Next, define
\begin{gather} 
    t_c'(x) \triangleq \Big\lbrace \argmax_{t\geq 0}\psi_h(t; x, u') \Big\rbrace \,, \label{eq:tc_prime} \\
    H'(x) \triangleq \psi_h(t_{c,0}'; x, u'), \; \forall t_{c,0}'\in t_c'(x) \,, \label{eq:H_prime} \\
	S_{H'} \triangleq \{x\in\reals \mid H'(x) \leq 0\} \,. \label{eq:SH'_set}
\end{gather}
\begin{remark}
	The functions $t_c'(x)$ and $H'(x)$ always exist because $\psi_h(t; x, u')$ is a polynomial with strictly negative highest coefficient $-a_{max}$.
\end{remark}
By definition, $H'(x) \geq h(x), \forall x \in \reals^n$, so $S_{H'} \subseteq S$ as well.
Lastly, we define the set $U'$,
\begin{equation}
	U' \triangleq \big\lbrace u \in U \;\big|\; \exists x \in S: u \in \boldsymbol{\mu}(x) \big\rbrace \subseteq U \,,
\end{equation}
which represents the set of control inputs such that $h^{(r)}=-a_{max}$. For brevity, in this section let $\psi_h(t) = \psi_h(t; x, u')$. 
We now state the second main result of this paper.

\begin{theorem} \label{thm:poly}
	The function $H'$ in \eqref{eq:H_prime} is a ZCBF on the set $S_{H'}$ in \eqref{eq:SH'_set} for the control set $U'$ (or $U$), provided $S_{H'} \neq \emptyset$.

	\begin{proof}
		As in \ref{thm:most_general}, $t_c'(x)$ is not necessarily unique, so we consider two cases depending on the elements of $t_c'(x)$. If $0\in t_c'(x)$, then $\psi_h(0) = h(x)$ and the ZCBF condition is satisfied by the same argument as \ref{thm:most_general}. 
		
		If $t_{c,0}'\in t_c'(x)$ for $t_{c,0}' > 0$, then $t_{c,0}'$ is a maximizer on an open interval so $\frac{\partial \psi_h}{\partial t}\big|_{t_{c,0}'} = 0$, where
		\begin{equation}
			\hspace{-4pt} 0= \left.\frac{\partial \psi_h}{\partial t}\right|_{t_{c,0}'} = \sum_{i=0}^{r-2} \frac{(t_{c,0}')^i}{i!}h^{(i+1)}(x) - \frac{(t_{c,0}')^{r-1}}{(r-1)!}a_{max} \,. \hspace{-2pt} \label{eq:step}
		\end{equation}
		Then one has
		\begin{align}
			&L_f \psi_h (t_{c,0}') + L_g \psi_h (t_{c,0}') u \overset{\eqref{eq:phi_poly}}{=} \sum_{i=0}^{r-2}\frac{1}{i!} h^{(i+1)} (x) (t_{c,0}')^i \nonumber \\ & \;\;\;\;\;\;\;\;\;\;\;\;\;\;\;+ \frac{1}{(r-1)!}\[ L_f^r h(x) + L_g L_f^{r-1} h(x) u\] (t_{c,0}')^{r-1}  \nonumber \\
			& \;\;\; \overset{\eqref{eq:step}}{=} \frac{(t_{c,0}')^{r-1}}{(r-1)!}\[a_{max} +  L_f^r h(x) + L_g L_f^{r-1} h(x) u\] \,. \label{eq:H'_dot} 
		\end{align}
		By definition of $a_{max}$, the right hand side of \eqref{eq:H'_dot} can be rendered nonpositive by a $u\in U'\subseteq U$ independent of $t_{c,0}'$, so by the same argument as \ref{thm:most_general}, $H'$ satisfies the definition of a ZCBF.
	\end{proof}
\end{theorem}

Similar to \ref{thm:most_general}, \ref{thm:poly} provides a guarantee of at least one safe trajectory, but $t_c'(x)$ is computed via polynomial root-finding rather than ODE propagation, making \eqref{eq:H'_dot} easier to compute than \eqref{eq:Hdot_explicit}-\eqref{eq:the_ode} when implementing condition \eqref{eq:cbf_condition}. Note the polynomial degree depends only on the relative-degree of $h$ (usually $r\leq 4$ \cite{Integrator_CBFs}), and not on the state dimension.
Similar to $\dot{H}$, note that if $t_c'(x)$ ever has multiple elements, then $\dot{H}'$ is given by \eqref{eq:Hdot_explicit} with $t_c'$ in place of $t_c$ and where $\frac{\partial \psi_h(t_{c,0})}{\partial x}$ is easily derived from \eqref{eq:phi_poly}.
See also \cite{Integrator_CBFs} for some restrictions on what $t_c'$ can be when $r=3, 4$. 
When $r=2$, $t_{c}'$ always contains only one element, which allows us to construct the following explicit form for $H'$.

\begin{example}
	In the case where $h\in\highdegree^2$, $H'$ takes the form:
	\begin{equation}
		H'(x) = \begin{cases}
			h(x) & \dot{h}(x) < 0 \\
			h(x) + \frac{\dot{h}(x)^2}{2a_{max}} & \dot{h}(x) \geq 0
		\end{cases} \,. \label{eq:the_example}
	\end{equation}	
\end{example}

\section{Case Study for Spacecraft Application} \label{sec:case_study}

In this section, we present two use cases for the CBFs in Section~\ref{sec:main_result} for a spacecraft in weak gravity, in which the spacecraft must navigate around an object under observation using control inputs calculated online via a ZCBF. 

In each case, the spacecraft state is $x = [r,\; v] \in \reals^6$ with
\begin{equation}
	\dot{x} = \begin{bmatrix} \dot{r} \\ \dot{v} \end{bmatrix} = \begin{bmatrix} v \\ f_\mu(r) + u \end{bmatrix} \,, \label{eq:dynamics}
\end{equation}
where $f_\mu$ is the local gravitational force. We construct a preplanned path $r_p(s), s\in[0,1],$ of desired observations which circumnavigate the observed object. This path is on the surface of the object, which is outside the safe set $S$, so the spacecraft must get close to $r_p$ while staying within the safe set. 
The spacecraft is driven to track the target using the Control Lyapunov Function (CLF) \cite{Exponential_CBF}: \,.
\begin{equation}
    V(x) = \frac{1}{2}||r-r_p||^2 + \frac{1}{2}k_2 ||v-k_1(r-r_p)||^2 \,.
\end{equation}

The spacecraft control input is then calculated as:
\begin{align}
	u(x) = &\argmin_{u\in U, \delta \in \reals} u\transpose u + J \delta^2 \;\;\;\;\textrm{such that} \label{eq:the_qp}\\
	&\;\;\;\;\; L_f H(x) + L_g H(x) u \leq \alpha(-H(x)) \nonumber \\
	&\;\;\;\;\; L_f V(x) + L_g V(x) u + \delta \leq - k_3 V(x) \nonumber
\end{align}
where $\delta$ is a slack variable for the CLF to ensure feasibility, $J$ is a constant slack penalty, $k_1,k_2,k_3$ are constants, and we assume $u$ in \eqref{eq:the_qp} is locally Lipschitz continuous. For simplicity, let $\alpha(\lambda) = \lambda$ in 
\eqref{eq:the_qp}.

\subsection{Case 1: Spherical Target}

\begin{figure}
	\centering
	\includegraphics[trim={1.9in, 0.2in, 1.4in, 1in},clip,width=0.8\columnwidth]{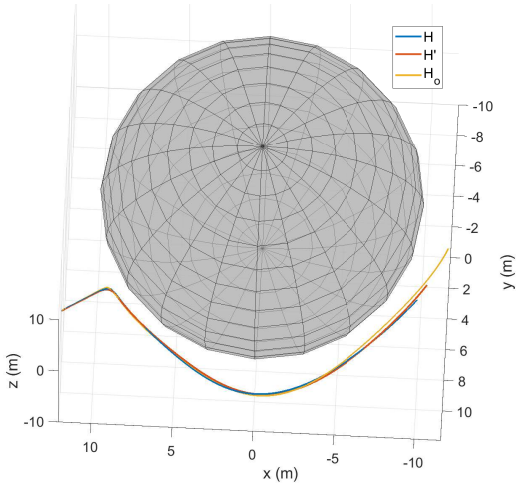}
	\caption{The paths around the spherical obstacle under the three ZCBFs considered ($H(x), H'(x), H_o(x)$)\vspace{-14pt}}
	\label{fig:sphere_trajectory}
\end{figure}

We first consider a spherical object of radius $\rho_t$ with fixed position $r_s$, and negligible gravity $f_\mu \equiv 0$, so \eqref{eq:dynamics} reduces to a double integrator. Maintaining a distance of $\rho_s$ from the object is equivalent to maintaining the relative-degree 2 constraint function $h_a (x) \leq 0$:
\begin{equation}
    h_a(x) = \rho_a - ||r - r_s|| \,,
\end{equation}
%
where $\rho_a = \rho_s + \rho_t > 0$. Note that $h_a$ is defined specifically so that $\ddot{h}_a$ represents physical acceleration, since if fuel is consumed slowly enough, then the spacecraft peak acceleration in any direction is a known constant, so $a_{max}$ is easy to compute. Suppose the spacecraft has six identical and orthogonal thrusters, so $U = \{u\in\reals^3 : ||u||_\infty \leq u_{max}\}$.


We then construct a ZCBF $H$ as in Section~\ref{sec:general} using the assumed control law $u^*_\textrm{ball}$ in \eqref{eq:opt_control_law}.
Specifically, to compute $H(x)$, we numerically propagate the dynamics \eqref{eq:dynamics} for some amount of time $T_\textrm{prop}$, yielding an array of states $\{\psi_x(t_k; x,u^*_\textrm{ball})\}_{k=1}^N$, and an array of $h$ values $\{\psi_h(t_k; x,u^*_\textrm{ball})\}_{k=1}^N$. For this particular combination of dynamics \eqref{eq:dynamics} and control law \eqref{eq:opt_control_law}, there is always a unique maximizer time $t_{c,0}$. Moreover, $t_{c,0}$ can be upper bounded, and we use this bound to choose the propagation time $T_{prop}$. Unfortunately, depending on the integration method, the array $\{\psi_h(t_k; x,u^*_\textrm{ball})\}_{k=1}^N$ may not include the true maximizer, since this array is only a sampling of a continuous curve. Thus, we select the three highest points $\{\psi_h(t_{k_l}; x,u^*_\textrm{ball})\}_{l=1}^{3}$ and fit a quadratic curve $q(\tau)=a \tau^2 + b\tau + c$ to these three points (where $\tau$ is time since $x$, i.e. $q(0) = h(x)$). Assuming the maximum is not simply $H(x) = h(x)$ (i.e. $t_{c,0}=0$), we then choose $H(x) = \max_{\tau\in\reals} q(\tau) = c - \frac{b^2}{4a}$, and denote $t_q = \argmax_{\tau\in\reals}q(\tau) = \frac{-b}{2a}$. From here, the gradient of $H$ can be computed as in \ref{remark:calculate_dphi} using $t_{c,0} = t_q$. Alternatively, for these dynamics specifically, it holds that $\theta(t_{c,0}) = \begin{bmatrix} I_{3\times3} & t_{c,0} I_{3\times3} \\ 0_{3\times3} & I_{3\times3} \end{bmatrix}$ in \eqref{eq:grad_psi}. For more details, see the simulation code below. Note, for more complex problems where $t_c(x)$ may have more than one element, if there are multiple local maximizers $t_k$ such that the values $\psi_h(t_k; x,u^*)$ are close to each other, then we recommend constructing $q(\tau)$ for each local maximizer.

Next, for the same system, we construct $H'$ as in Section~\ref{sec:poly}. It follows that $H'$ is as given in \eqref{eq:the_example} with $h_a$ in place of $h$ and $a_{max} = u_{max}$. 
This leads to $U' = \{u\in\reals^3:||u|| \leq u_{max}\}\subset U$. 
Finally, for comparison, consider the ZCBF 
\begin{equation}
	H_o(x) = \left(\arctan(\dot{h}_a(x))  + \frac{\pi}{2} \right)h_a(x) \label{eq:ames_cbf} \,,
\end{equation}
derived using the rules in \cite{Grizzle_SeminalPaper}. As discussed in \cite{Grizzle_SeminalPaper}, and similar to \cite{Backstepping_CBF,Exponential_CBF,XU2018195}, this ZCBF is only valid over $U = \reals^3$.

\begin{figure}
	\centering
	\begin{minipage}{\columnwidth}
		\centering
		\includegraphics[width=\textwidth]{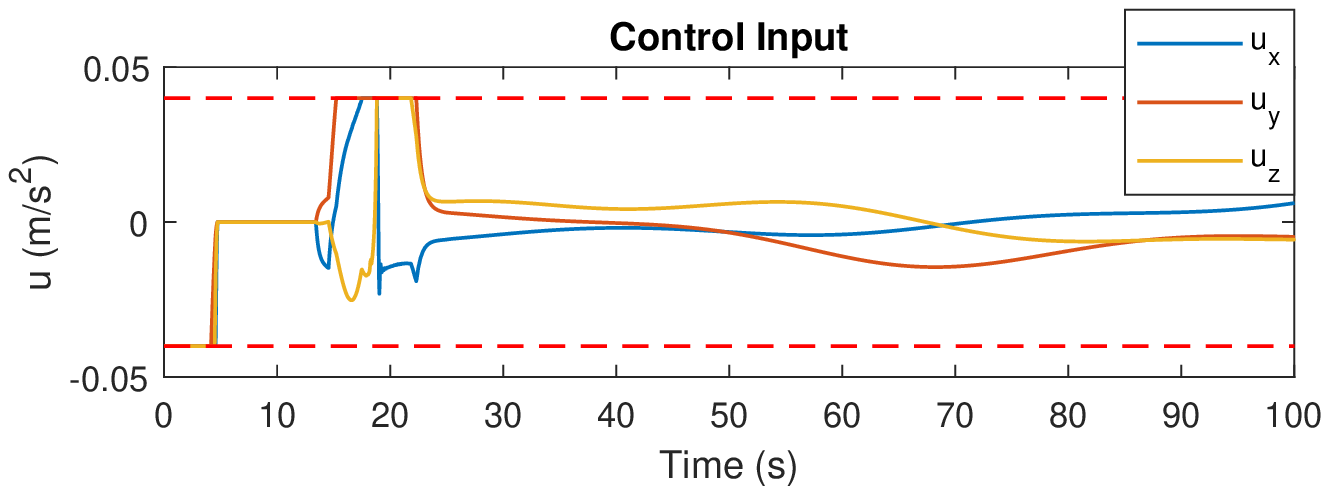}
		\captionof{figure}{The control input using $H(x)$ as the ZCBF}
		\label{fig:sphere_inf}
	\end{minipage}\vspace{6pt}
	\begin{minipage}{\columnwidth}
		\centering
		\includegraphics[width=\textwidth]{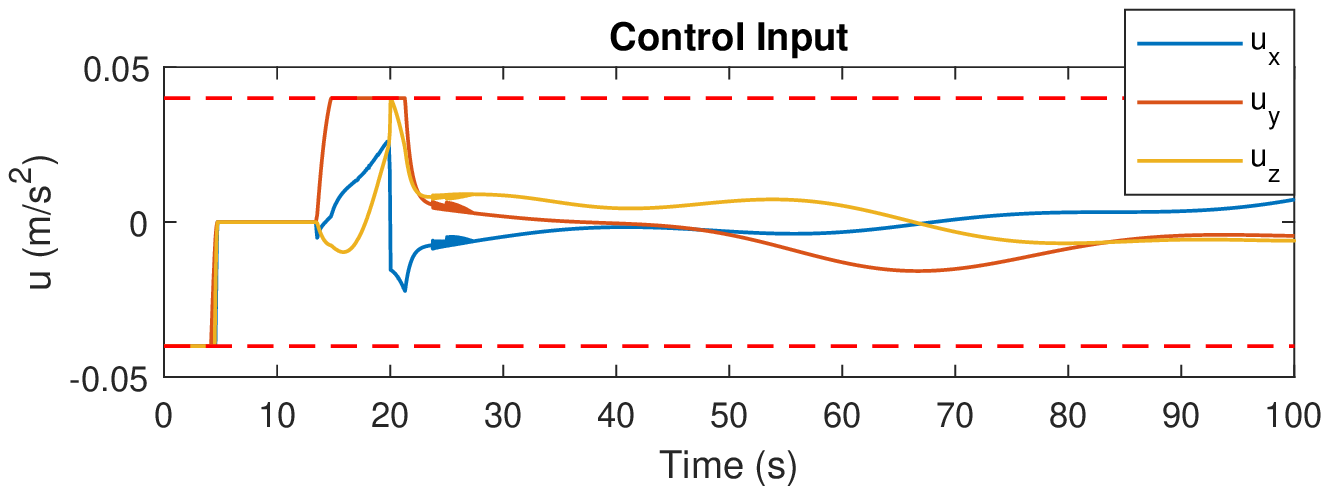}
		\captionof{figure}{The control input using $H'(x)$ as the ZCBF}
		\label{fig:sphere_two}
	\end{minipage}\vspace{6pt}
	\begin{minipage}{\columnwidth}
		\centering
		\includegraphics[width=\textwidth]{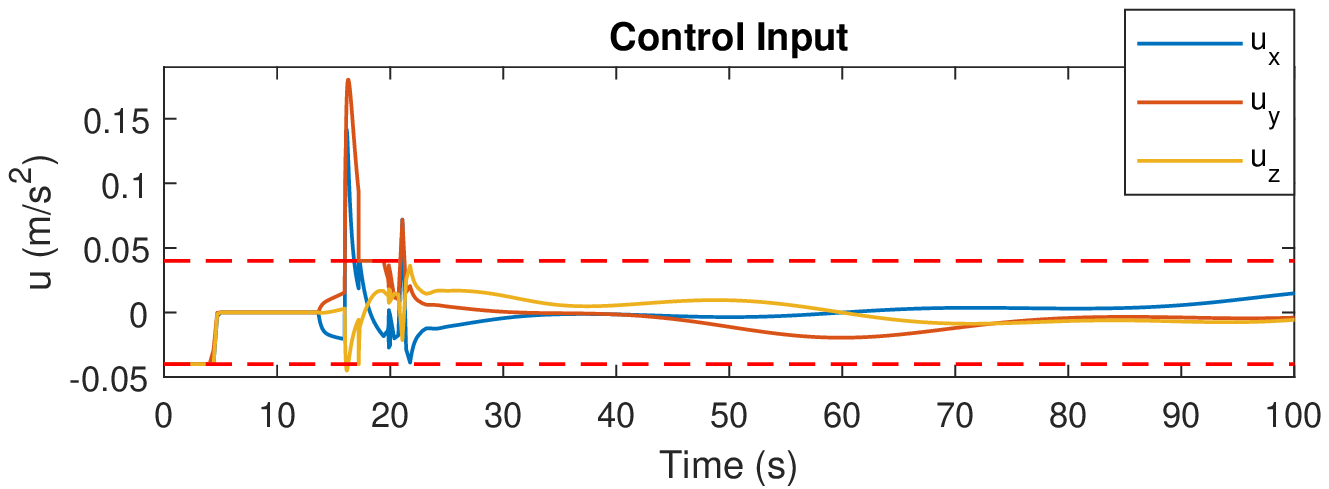}
		\captionof{figure}{The control input using $H_o(x)$ as the ZCBF, which necessitates using control inputs outside the prescribed bounds (dashed red lines) for \eqref{eq:the_qp} to have a solution \vspace{-8pt}}
		\label{fig:sphere_compare}
	\end{minipage}
\end{figure}

We then simulated trajectories under all three ZCBFs ($H$, $H'$, $H_o$) in MATLAB\footnote{All simulation code can be found at \url{https://github.com/jbreeden-um/phd-code/tree/main/2021}}. The three simulated paths are shown in Fig.~\ref{fig:sphere_trajectory}, and the corresponding control inputs are shown in Figs.~\ref{fig:sphere_inf}-\ref{fig:sphere_compare}. All three trajectories in Fig.~\ref{fig:sphere_trajectory} remained safe and generally followed similar paths around the obstacle. However, the control input only stayed within the designated control set in Figs.~\ref{fig:sphere_inf}-\ref{fig:sphere_two}. For the trajectory under $H_o(x)$, the QP in \eqref{eq:the_qp} became infeasible at $t=14$ (see Fig.~\ref{fig:sphere_compare}), so we had to expand the control set $U$ to compute a control input satisfying the safety constraint in \eqref{eq:cbf_condition}. Thus, the proposed methods always yield trajectories that are safe in the presence of control input constraints, whereas earlier methods might not.

Additionally, note that the lines for $u_x, u_z$ for $t\in[15,25]$ are slightly closer to the horizontal axis in Fig.~\ref{fig:sphere_two} than in Fig.~\ref{fig:sphere_inf}. This occurred because $H'$ only makes use of $u\in U'\subset U$ even though the QP was calculated over the complete control set $U$.

\subsection{Case 2: Asteroid Target} \label{sec:asteroid_study}

\begin{figure}
	\centering
	\includegraphics[trim={0.9in, 0.2in, 1.2in, 0.5in},clip,width=\columnwidth]{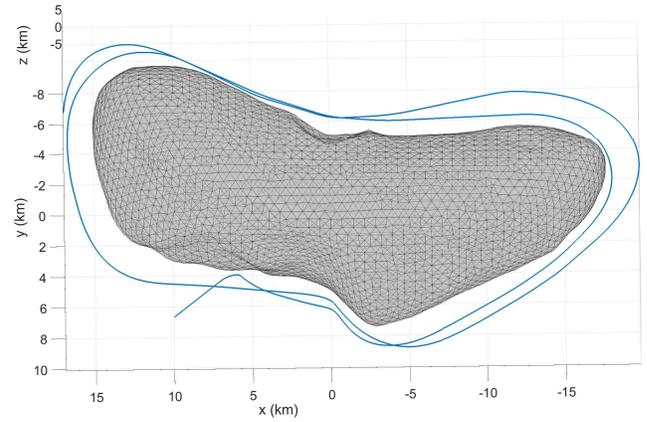}
	\caption{The trajectory around Eros (top view)\vspace{-6pt}}
	\label{fig:asteroid_trajectory}
\end{figure}

\begin{figure}
	\centering
	\includegraphics[width=\columnwidth]{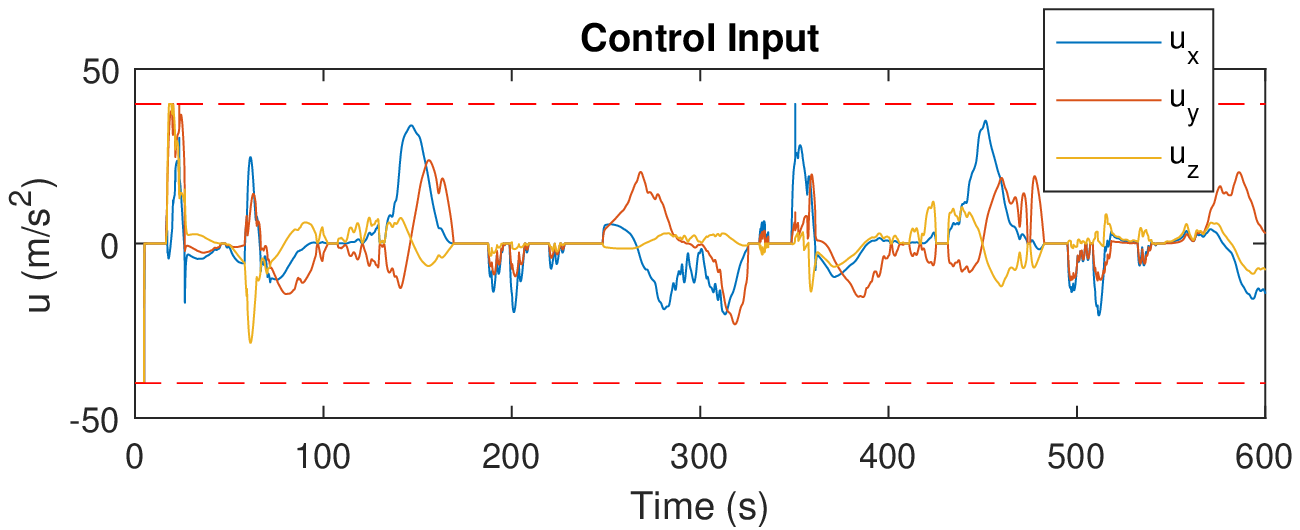}
	\caption{The control input over time for Fig.~\ref{fig:asteroid_trajectory}, which stays between the prescribed bounds (dashed red lines) \vspace{-8pt}}
	\label{fig:asteroid_control}
	\vspace{16pt}
	\includegraphics[width=3.3in]{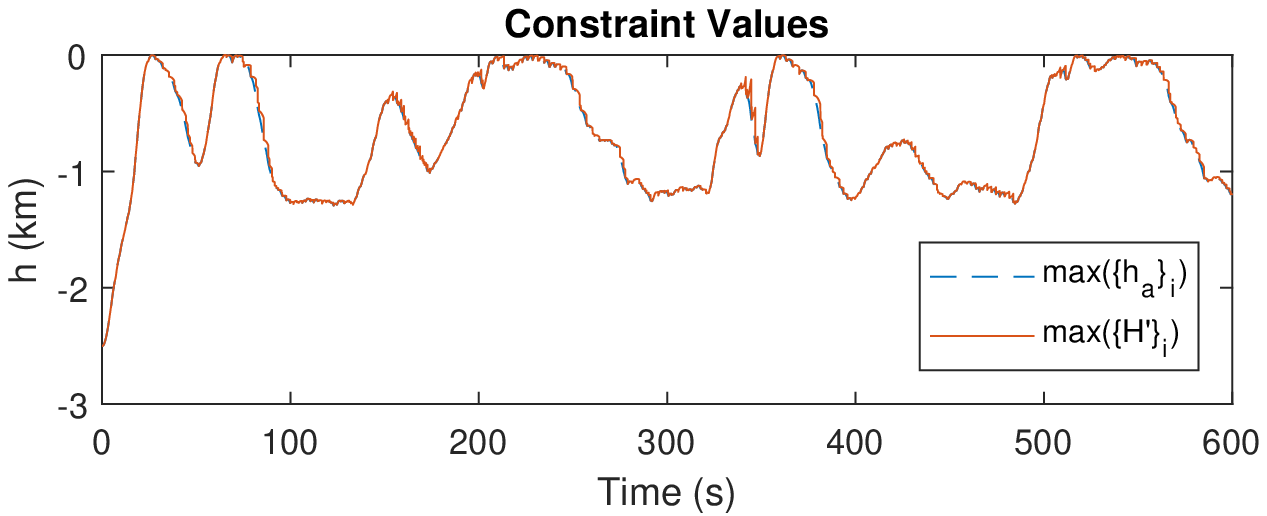}
	\caption{The value of the largest $h_a(x)$ and $H'(x)$ over time\vspace{-8pt}}
	\label{fig:asteroid_constraint}
\end{figure}

Next, we consider a spacecraft avoiding an object more complicated than a sphere, in this case an asteroid, but still in weak gravity, with gravity modeled by spherical harmonics \cite{Healy}. If the asteroid is convex, one could simply use $h_a(x)$ as the constraint function, setting $r_s$ as the instantaneous closest point. 
However, if the asteroid is nonconvex (as in this example), the spacecraft could obtain a large velocity with respect to a point other than the closest point, so that strategy is no longer sufficient. Instead, we consider a discrete point-cloud model of the entire asteroid. The simplest response is to construct a ZCBF for every point in the model, similar to \cite{Intersecting_CBFs}, though this could be computationally demanding. The result in \cite[Thm. 3]{NonSmooth_CBFs} allows us to to reduce the number of constraints to only those constraints that could be violated within some finite time horizon $\Delta t > 0$ to reduce complexity.
For this case study, we chose asteroid 433 Eros. The point cloud is a shape model of Eros with $N=7790$ plates \cite{shape_model}, and gravity $f_\mu$ modeled in \cite{gravity_model}. 
We construct a ZCBF with the form of $H'(x)$ for every point,
yielding $\lbrace H_i'(x) \rbrace_{i=1}^N$. 
The trajectory, control inputs, and CBF values for the asteroid simulation are shown in Figs.~\ref{fig:asteroid_trajectory}-\ref{fig:asteroid_constraint} and a video of the scenario can be found at \url{https://youtu.be/JKj3PUrYnEg}. As expected, the trajectory is always safe and $u$ satisfies the input constraints. Since the spacecraft generally moves tangentially to the asteroid, $H$ is only slightly larger than $h_a$ for most of the simulation.

\section{Conclusions} \label{sec:conclusions}

We have introduced two ZCBF formulations for high-relative-degree constraint functions. Both approaches are derived from control policies known to meet the input constraints, and thus are guaranteed safe under such constraints without further tuning. The two strategies were demonstrated safe under input constraints on a spacecraft obstacle-avoidance example, whereas most prior techniques would have failed to meet the input constraints \cite{Grizzle_SeminalPaper,Backstepping_CBF,Exponential_CBF,My_Paper} or required clever selection of bounding functions \cite{CBFs_ManyClassK}. The second strategy was further verified on an asteroid exploration example. 
Future work includes extending this formulation to guarantee feasibility of more complex objectives and unsafe sets, in particular, environments with strong central gravity, and comparing the efficiency of online approaches to optimal planned paths. 


\bibliographystyle{ieeetran}
\bibliography{sources}

\end{document}